\documentclass[12pt]{article}
\usepackage[french]{babel}
\usepackage{amssymb}
\usepackage{eufrak}
\usepackage{amsmath}
\usepackage{t1enc}
\usepackage{mathrsfs}
\usepackage[latin1]{inputenc}

\include{xy}
\xyoption{all}

\begin{document}
\newtheorem{theorem}{Theorem}[section]
\newtheorem{lemme}[theorem]{Lemma}
\newtheorem{il}[theorem]{Illustration}
\newtheorem{coro}[theorem]{Corollary}
\newtheorem{Prop}[theorem]{Proposition}
\newtheorem{defin}[theorem]{Definition}
\newtheorem{remark}[theorem]{Remark}
\newtheorem{remarks}[theorem]{Remarks}
\newtheorem{ex}[theorem]{Examples}
\newtheorem{hyp}[theorem]{Hypothesis}
\newtheorem{Not}[theorem]{Notations}

\begin{center}

Computation of the Galois groups\\
 occurring in  M. Papanikolas's study of   Carlitz
logarithms
\end{center}

\begin{center}

Charlotte Hardouin (IWR, Heidelberg)

\end{center}

\begin{center}
February 23, 2007

\end{center}

\section{Introduction}
In  this note, we give an alternative presentation of one of the ingredients occurring in M. Papanikolas's proof  of the algebraic
independence of Carlitz logarithms  
\cite{Papa}. More precisely, the main theorem of \cite{Papa} reduces the problem to the computation of the Galois group $G_X$ of a certain $t$-motive $X$, and we present an alternative proof of the computation of  $G_X$. The method is inspired from \cite{Hardouin},
 and  should apply to other situations, such as logarithms of Drinfeld elliptic modules, or values of $\zeta$-functions. We now recall the statement of Papanikolas's theorem, and the notations of his article, which will be kept throughout this note. \\
\subsection{Notations}
 Let $\mathbb{F}_q$ the field of $q$-elements, where $q$ is a prime
power of $p$. Let $k=\mathbb{F}_q(\theta)$, where $\theta$ is
transcendental over $\mathbb{F}_q$, and define an absolute valuation
$|.|_{\infty}$ at the infinite place of $k$  such that
$|\theta|_{\infty}=q.$ let $k_{\infty}$ be the $\infty$-adic
completion of $k$, let $\overline{k_{\infty}}$ be an  algebraic
closure, let $\mathbb{K}$ be the $\infty$-adic completion of
$\overline{k_{\infty}}$, $\mathbb{T}:= \mathbb{K}\lbrace t \rbrace$
is the ring of restricted power series and let $\overline{k}$ be the
algebraic
closure of $k$ in $\mathbb{K}$. For $f = \sum_i a_i t^i$ in $\mathbb{T}$, we set $f^{(-1)} = \sum_i a_i^q t^i$.\\

\begin{defin} [see \cite{Papa}] \label{theorem:rigid} 
We let   $\mathcal{T}$ be the category of $t$-motives in the sense of \cite{Papa}, 3.4.10. \end{defin}

We recall that $\mathcal{T}$ is a  strictly full Tannakian sub-category of the category $\mathcal{R}$ of  rigid analytically trivial pre-$t$-motives. Objects in $\mathcal{R}$ correspond to certain $\sigma$-difference equations over  $\overline{k}(t)$, and a fiber functor $\omega$ on $\mathcal{T}$ is provided by rigid analytic trivialization. In particular, 
$\mathcal{T}$ is a neutral tannakian category  over
$\mathbb{F}_q(t)$. We denote its identity object by $\bold{1}$, and for any $X$ in $\mathcal{T}$, we write $G_X = Aut^{\otimes}(\omega_{|<X>})$ for the Galois group of $X$ attached to the fiber functor $\omega$, see \cite{Papa}, 3.5.2, 4.4.1 and  5.4.10.  By \cite{Papa}, 5.2.12.b, this is a reduced affine group scheme over $\mathbb{F}_q(t)$.

\subsection{Exemples of $\sigma$-equations associated to objects of
$\mathcal{T}$}

\begin{enumerate}
\item \textbf{The Carlitz motive}\\

We define the Carlitz motive to be the pre-$t$-motive $\cal C$ 
whose
underlying $\overline{k}(t)$-vector space is $\overline{k}(t)$
itself and on which $\sigma$ acts by
$$ \sigma(f):=(t-\theta)f^{(-1)}, f \in {\cal C}$$
\begin{enumerate}
\item The Carlitz motive is rigid analytically trivial  and one of
its analytic trivialization is given by the function $1\over \Omega$ (see
\cite{Papa} 3.3.5).
\item The number $\tilde{\pi}= -\frac{1}{\Omega(\theta)}$ is the
\textit{Carlitz period}.
\item The  Galois group $G_{\cal C}$ of $\cal C$
is equal to $\mathbb{G}_m$.
\item Moreover, we have $End_{\mathcal{T}}(C,C)=\mathbb{F}_q(t)$ (cf. \cite{Papa} 3.5.3).
\end{enumerate}

\medskip

\item \textbf{The Carlitz logarithm motive}\\

\smallskip

Let $\alpha_i \in \overline{k}^*$ with $|\alpha_i|_{\infty}<
|\theta|_{\infty}^{q/(q-1)}$. Set :
$$\Phi(\alpha_i) :=\left( \begin{array}{cc}
(t-\theta) & 0 \\
\alpha_i^{(-1)} (t-\theta) & 1 \end{array} \right).$$
$\Phi(\alpha_i)$ defines a pre-$t$-motive $X(\alpha_i)$, which is an
extension in  the category $\mathcal{T}$ of $\bold{1}$ by the Carlitz motive $\cal C$
$$\xymatrix{
0 \ar[r] & {\cal C} \ar[r] & X(\alpha_i) \ar[r] & \bold{1}  \ar[r] & 0} .$$
 Indeed, the
pre-t-motive $X(\alpha_i)$ is rigid analytically trivial (see
\cite{Papa} prop. 7.1.3) and its trivialization is given by:
$$\Psi(\alpha_i) :=\left( \begin{array}{cc}
\Omega & 0 \\
\Omega L_{\alpha_i} & 1 \end{array} \right),$$ where the function
$L_{\alpha_i}$ is defined as in \cite{Papa}, 7.1.1 : this is an element of $\mathbb{T}$
 satisfying the functional equation :
 $$ \sigma(L_{\alpha_i})=\alpha_i^{(-1)} +
 \frac{L_{\alpha_i}}{(t-\theta)},$$ whose value  at $t = \theta$ is equal to the Carlitz
logarithm  $Log_{\cal C}(\alpha_i)$ of $\alpha_i$.

\item \textbf{The multiple Carlitz logarithm motive}\\

Let $\alpha_1,...,\alpha_r \in \overline{k}^*$ with
$|\alpha_i|_{\infty}< |\theta|_{\infty}^{q/(q-1)}$. Set :

$$ \Phi(\alpha_1,...,\alpha_r):= \left( \begin{array}{cccc}
t-\theta & 0 & \cdots & 0 \\
\alpha_1^{(-1)} (t-\theta)& 1& \cdots & 0 \\
\vdots & \vdots & \ddots & \vdots \\
 \alpha_r^{(-1)} (t-\theta)& 0 & \cdots &1 \end{array} \right).$$
$\Phi(\alpha_1,...,\alpha_r)$ defines a pre-$t$-motive $X(\alpha_1,
...,\alpha_r)$ which is an extension of $\bold{1}^r$ by the Carlitz
motive $\cal C$
:
$$\xymatrix{
0 \ar[r] & {\cal C} \ar[r] & X(\alpha_1,...,\alpha_r) \ar[r] & \bold{1}^r
\ar[r] &0 }.$$ The pre-t-motive $X(\alpha_1,...,\alpha_r)$ is rigid
analytically trivial (see \cite{Papa} prop. 7.1.3) and its
trivialization is given by:
$$ \Psi(\alpha_1,...,\alpha_r):= \left( \begin{array}{cccc}
\Omega & 0 & \cdots & 0 \\
\Omega L_{\alpha_1} 1& \cdots & 0 \\
\vdots & \vdots & \ddots & \vdots \\
 \Omega L_{\alpha_r}& 0 & \cdots &1 \end{array} \right).$$

\end{enumerate}

As in \cite{Hardouin}, proof of Cor. 2.2,   we have:
\begin{lemme}\label{lemme:cat}
The tannakian category generated by $X(\alpha_1,...,\alpha_r)$ in
$\mathcal{T}$ is equal to the Tannakian category generated by the
motive $\bigoplus_{i=1}^r X(\alpha_i)$.
\end{lemme}

\subsection{Papanikolas's theorems on algebraic independence.}

\begin{theorem}[Theorem 7.4.2 in \cite{Papa}]\label{theorem:ilc}
Let $\lambda_1,...,\lambda_r \in \mathbb{K}$ satisfy
$exp_C(\lambda_i) \in \overline{k}$ for $i=1,...,r$. If
$\lambda_1,...,\lambda_r$ are linearly independent over $k$, then
they are algebraically independent over $\overline{k}$.
\end{theorem}
Since the period ${\tilde {\pi}}$  satisfies $exp_{\cal C}({\tilde \pi}) = 0$, we can rephrase Theorem \ref{theorem:ilc} as
follow :
Let $\lambda_1,...,\lambda_r \in \mathbb{K}$ satisfy
$exp_C(\lambda_i) \in \overline{k}$ for $i=1,...,r$. If
$\lambda_1,...,\lambda_r,\tilde{\pi} $ are linearly independent over
$k$, then they are algebraically independent over $\overline{k}$.

Because the indetermination of the Carlitz logarithm is given by $k$-multiples of $\tilde \pi$ (cf. \cite{Papa}, 7.4.1), this is in turn equivalent to 

\begin{theorem} \label{theorem:ilc2}
Let $\alpha_1,...,\alpha_r \in \overline{k}^*$ with
$|\alpha_i|_{\infty}< |\theta|_{\infty}^{q/(q-1)}$. Assume that $\tilde{\pi},  log_{\cal C}(\alpha_1),..., log_{\cal C}(\alpha_r)$
 are linearly
independent over $k$. Then they are algebraically independent over
$\overline{k}$.
\end{theorem}

Now, $\tilde{\pi} = -\frac{1}{\Omega(\theta)},  log_{\cal C}(\alpha_1) = L_{\alpha_1}(\theta),..., log_{\cal C}(\alpha_r) = L_{\alpha_r}(\theta)$. Combining the main Theorem 1.1.7 of his article together with a previous transcendence criterion (Theorem 6.1.1), Papanikolas reduces the proof of Theorem \ref{theorem:ilc2} to showing:

\begin{theorem}[Theorem 7.3.2.c in \cite{Papa}] \label{theorem:ilc3}
Let $\alpha_1,...,\alpha_r \in \overline{k}^*$ with
$|\alpha_i|_{\infty}< |\theta|_{\infty}^{q/(q-1)}$. Assume that $\tilde{\pi},  log_{\cal C}(\alpha_1),..., log_{\cal C}(\alpha_r)$
 are linearly
independent over $k$. Then the dimension of the Galois group $G_X$ of the $t$-motive $X = X(\alpha_1,
...,\alpha_r)$ is equal to $r + 1$.
\end{theorem}

\subsection{Sketch of the proof of Theorem \ref{theorem:ilc3} }

Following \cite{Papa}, we will work in the framework of the Tannakian
category $\mathcal{T}$ of $t$-motives, cf. Definition \ref{theorem:rigid}. As just recalled, the method of M. Papanikolas for proving Theorem
\ref{theorem:ilc} is to compute the Galois group $G_X$ of the $t$-motive $X$. This is the content of Theorem 7.3.2 of \cite{Papa}, where $G_X$ is denoted by $\Gamma_X$. Note, however, that the paragraph following $(7.2.4.1)$ needs some clarification,  since $\Gamma_X$ is not a linear subspace.
  In this
note, we will give a tannakian version of the computation of $G_X$, which while settling this point,  actually simplifies  the proof of \cite{Papa}, and points towards further generalizations of  Theorem \ref{theorem:ilc}.\\

So, we have to compute the dimension of the Galois group attached to
the motive $X =   X(\alpha_1,
...,\alpha_r)$. To this purpose, we deduce from  Lemma
\ref{lemme:cat} that the Galois group of $X$ is equal to the Galois
group $G$ of $\bigoplus_{i=1}^r X(\alpha_i)$. As in \cite{Papa}, 7.2.2, we  see that the quotient
of $G$ by its unipotent radical is isomorphic to the Galois group of
the Carlitz motive $\cal C$, i.e to $\mathbb{G}_m$. Therefore, it remains
to compute the dimension of the unipotent radical of $G$, that
is the unipotent radical of the Galois
group of a sum of extensions of $\bold{1}$ by the Carlitz motive.\\

To compute the latter dimension, we will use the theorems of Section $2$ below, which reduce the problem  to a question
of linear algebra; this section combines the arguments of \cite{Hardouin} with Papanikolas's  crucial observation that the unipotent radical is a {\it vectorial}  group, see \cite{Papa}, 7.2.3.  
Finally, Section 3  completes the proof of Theorem \ref{theorem:ilc3}, along the lines of  \cite{Papa}, bottom of p. 50.

\section{Computation of Galois groups in Tannakian categories in
characteristic $p$}  

Let $p$ be a prime number. Let $(\bold{T}, \omega)$ be a neutral Tannakian
category over a field $C$
 of characteristic $p$. Let $\bold{1}$
denotes the unit object of $\bold{T}$, so that $C=End(\bold{1})$ and $\omega : \bold{T} \mapsto
Vect_C$. In the application to \cite{Papa}, $\bold{T}= \mathcal{T}$ and
$C=\mathbb{F}_q(t)$, where $q$ is a power of $p$ and $t$ is
transcendental over $\mathbb{F}_q$ .
\\

For any object $\mathcal{X}$ in  $\bold{T}$, we denote by   $G_\mathcal{X}$
  the linear algebraic group
scheme $Aut^{\otimes}(\omega|_{<\mathcal{X}>})$ over $C$. Furthermore, we identify  $C$-vector spaces such as $\omega(\mathcal{X})$ to {\it  vectorial groups}  over $C$.   

\begin{theorem}\label{theorem:Ber2}

Let $\mathcal{Y}$ be an object of $\bold{T}$,  and let $\mathcal{U}$ be an extension of $\bold{1}$ by $\mathcal{Y}$. Assume that $G_\mathcal{U}$ is reduced,
  that $G_\mathcal{Y} = \bold{G}_m$, and that the action of ${\bold{G}_m}$ on $\omega(\mathcal{Y})$ is given by its canonical character.
 Then
the unipotent radical of the Galois group $G_\mathcal{U}$ is equal
to $\omega(\mathcal{V})$ where $\mathcal{V}$ is the smallest
sub-object of $\mathcal{Y}$ such that $\mathcal{U}/\mathcal{V}$ is a
trivial extension of $\bold{1}$ by $\mathcal{Y}/\mathcal{V}$.
\end{theorem}

\textbf{Proof}\\
First of all, we remark that every $\bold{G}_m$-module of
finite dimension over $C$ is completely reducible (see \cite{Jantzen}
p.35). By Tannaka theorem (see \cite{Dlct}) , there is an
equivalence of category between $<\mathcal{Y}>$ and the category
$Rep_{G_\mathcal{Y}}$ of $G_\mathcal{Y}$-modules of finite dimension over $C$. Then, it
is clear that $\mathcal{Y}$ is a completely reducible object in
$\bold{T}$.

\medskip

\textit{Existence of the smallest sub-object}\\

Let us denote by
$\bold{V}$ the set of sub-objects $\mathcal{W}$ of $\mathcal{Y}$
such that $\mathcal{U}/ \mathcal{W}$ is a trivial extension of
$\bold{1}$ by $\mathcal{Y}/\mathcal{W}$. It is enough to prove that
if $\mathcal{V}_1$ and $\mathcal{V}_2$ are in $\bold{V}$, their
intersection $\mathcal{W}$ lies in $\bold{V}$.\\

Because $\mathcal{Y}$ is completely reducible, there exist three
sub-objects $\mathcal{V'}$, $\mathcal{W}_1'$, $\mathcal{W}_2'$ of
$\mathcal{Y}$ such that :
\begin{enumerate}
\item $\mathcal{V}_1=\mathcal{W} \oplus \mathcal{W}_1'$, $ \mathcal{V}_2=\mathcal{W} \oplus \mathcal{W}_2'$.
\item $\mathcal{Y} = \mathcal{V}_1 \oplus \mathcal{W}_2'\oplus \mathcal{V}' = \mathcal{V}_2 \oplus \mathcal{W}_1' \oplus \mathcal{V}'=\mathcal{W} \oplus \mathcal{W}_2' \oplus \mathcal{W}_1' \oplus \mathcal{V}$
\end{enumerate}
We have : $$Ext^1(\bold{1}, \mathcal{Y})\simeq Ext^1(\bold{1},
\mathcal{V}_1) \times Ext^1(\bold{1},\mathcal{W}_2'\oplus
\mathcal{V'}) \ \mbox{ et} \  Ext^1(\bold{1}, \mathcal{Y})\simeq
Ext^1(\bold{1}, \mathcal{V}_2) \times
Ext^1(\bold{1},\mathcal{W}_1'\oplus \mathcal{V'}) .$$ Because
$\mathcal{V}_1$ and $\mathcal{V}_2$ are in  $\bold{V}$, the
projection of  $\mathcal{U}$ is  trivial on  $Ext^1(\bold{1},
\mathcal{W}_2'\oplus \mathcal{V'})$ and on $
Ext^1(\bold{1},\mathcal{W}_1'\oplus \mathcal{V'} )$. Then the
projection of $ \mathcal{U}$ is also trivial on
$Ext^1(\bold{1},\mathcal{W}_2'\oplus \mathcal{W}_1'\oplus
\mathcal{V'})$ and thus $\mathcal{W}$ is in $\bold{V}$.\\

\textit{ Computation of the unipotent radical $R_u$ of the Galois
group $G_\mathcal{U}$ of $\mathcal{U}$}\\

By assumption, $\mathcal{U}$  lies in an exact sequence:
$$\xymatrix{
0 \ar[r] &\mathcal{Y} \ar[r]^{i} & \mathcal{U} \ar[r]^{p} & \bold{1}
\ar[r] & 0}.$$ Let $R$ be a $C$-algebra. Since the
categories $<\mathcal{U}>$ and $Rep_{G_\mathcal{U}}$ are equivalent,
$\omega(\mathcal{U})\otimes R$ is an extension of the unit
representation $1_R$ par $\omega(\mathcal{Y})\otimes R$ in the
category $Rep_{G_\mathcal{U}(R)}$ of $G_\mathcal{U}(R)$-modules of finite rank over
$R$. Consider the exact sequence
of  free $R$-modules : \\
$$\xymatrix{
0 \ar[r] &\omega(\mathcal{Y})\otimes R   \ar[r]^{\omega(i)^R} &
\omega(\mathcal{U})\otimes R & \ar[r]^{\omega(p)^R} & R
\ar@{.>}@/^/[l]^{s^R}
\ar[r] & 0} , $$ 
fix a section $s$ of the underlying exact sequence of $C$-vector spaces, and put  $f^R=s^R(1) \ \in \omega(\mathcal{U}) \otimes R$, where $s^R = s \otimes 1$.\\
Let us consider the morphism  of  $C$-schemes $ \zeta_{\omega(\cal
U)}^R : G_\mathcal{U}(R) \rightarrow \omega(\mathcal{Y})\otimes R$
 defined by
the relation :

$$ \forall \sigma \in G_\mathcal{U}(R),  \zeta_{\omega(\cal
  U)}^R(\sigma) = (\sigma-1)f^R.$$ This defines
 a morphism of schemes $\zeta_{\omega(\cal
  U)}$
 over $C$ from
$G_\mathcal{U}$ with value in the $C$-vector space
$\omega(\mathcal{Y})$, whose restriction to $R_{u}$ is an immersion
of algebraic group-schemes over $C$ from $R_u$ to the
$C$-vectorial group  $\omega(\mathcal{Y})$.
Since $G_\mathcal{U}$ is reduced, its scheme theoretic image  is again reduced,
and we have:

\medskip
\begin{lemme} [see \cite{Hardouin}, 2.8 and \cite{Papa}, 7.2.3] \label{lemme:norm}
The image $W$ of $R_u$ under $\zeta_{\omega(\cal U)}$ is a
$C$-vectorial subgroup of the $C$ vectorial group  
$\omega(\mathcal{Y})$.

\end{lemme}
 
\textbf{Proof}\\
Since $W$ is reduced, it suffices to check this on points in the algebraic closure of $C$.
For all $\sigma_1  \in G_\mathcal{Y}$ and $\sigma_2 \in R_u$, we
have
$$\zeta_{\omega(\cal U)}( \sigma_1 \sigma_2
{\sigma_1}^{-1})=\sigma_1(\zeta_{\omega(\cal U)}(\sigma_2)).$$

Indeed, we  have : \begin{equation} \label{eqn:co2} \sigma_1
\zeta_{\omega(\cal U)}({\sigma_1}^{-1})= (1 - \sigma_1) f =
-\zeta_{\omega(\cal U)}(\sigma_1),\end{equation}\\
and
$$\zeta_{\omega(\cal U)}( \sigma_1 \sigma_2
{\sigma_1}^{-1})= \sigma_1(\zeta_{\omega(\cal
U)}(\sigma_2{\sigma_1}^{-1})) + \zeta_{\omega(\cal U)}(\sigma_1)
=\sigma_1(\sigma_2(\zeta_{\omega(\cal
U)}({\sigma_1}^{-1}))+\zeta_{\omega(\cal
U)}(\sigma_2))+\zeta_{\omega(\cal U)}(\sigma_1).$$ From
(\ref{eqn:co2}), we deduce that :
$\sigma_1(\sigma_2(\zeta_{\omega(\cal
U)}({\sigma_1}^{-1})))=-\sigma_1 \sigma_2
{\sigma_1}^{-1}(\zeta_{\omega(\cal U)}(\sigma_1))$. But $\sigma_1
\sigma_2 {\sigma_1}^{-1}$ is an element of $R_u$ and
$\zeta_{\omega(\cal
U)}(\sigma_1)$ lies in $\omega(\mathcal{Y})$. 
 Then, $\sigma_1(\sigma_2(\zeta_{\omega(\cal U)}({\sigma_1}^{-1})))=
-\zeta_{\omega(\cal U)}(\sigma_1)$. Therefore
$\sigma_1(\zeta_{\omega(\cal U)}(\sigma_2))=\zeta_{\omega(\cal U)}( \sigma_1 \sigma_2 {\sigma_1}^{-1})$ belongs to  $W$.\\

In other words, $W$ is an algebraic subgroup over $C$ of $\omega(\mathcal{Y})$ which is stable under the action of $G_\mathcal{Y}$. Now,
$G_\mathcal{Y}=\mathbb{G}_m$ and the hypothesis that $\omega(\mathcal{Y})$ is an {\it isotypic} representation of $\mathbb{G}_m$ implies that W is  a $C$-vectorial
subgroup of the $C$-vectorial group $\omega(\mathcal{Y})$.
\medskip

 \begin{lemme} [see \cite {Hardouin}, 2.9] \label{lemme:minimal} The image under  $\omega$  of the smallest sub-object of $\bold{V}$
is equal to   $W$.
\end{lemme}

\textbf{Proof}\\
Let us denote by
  $\mathcal{V}$ the minimal object of $\bold{V}$, and by $V$ its image under $\omega$.
  Then, $G_\mathcal{U}$
acts on  $\omega(\mathcal{U}/\mathcal{V})$
through $G_\mathcal{Y}$ (because $\mathcal{U}/\mathcal{V}$ is a trivial extension of  $ \bold{1}$ by a quotient of $\mathcal{Y}$ in the category $\bold{T}$).
  Thus the projection of $f^C=s(1)$
  in $ \omega(\mathcal{U})/V$  is invariant under the action of  $R_u$, 
  and  the orbit $\lbrace \sigma f^C -f^C ; \ \sigma \in R_u \rbrace$
 lies in $V$.
   Therefore $\zeta_{\omega(\cal U)} (R_u):= W \subset V$.\\
   
Conversely,  the image   $W$ of $R_u$ under $\zeta_{\omega(\cal
U)}$ is, by Lemma \ref{lemme:norm}, a $C$-vector-space stable
under the action of $G_\mathcal{Y}$ in $\omega(\mathcal{Y})$.  Then, by
equivalence of category, there exists a sub-object $\mathcal{W}$ of $\mathcal{Y}$ in 
$\bold{T}$ such that $\omega( \mathcal{W})=W$. Let us show that $\mathcal{W} $ is an element
of $\bold{V}$. 
Since $W$ is the image of $R_u$, $G_{\mathcal{U}}$ acts on  $\omega(\mathcal{U}) / W$ through its quotient $G_{\mathcal{U}}/R_u = G_{\mathcal{Y}}$.
Therefore,  $\omega(\mathcal{U})/W(C)$ is  an  extension of $C$ by $\omega(\mathcal{Y})/W(C)$ in
the category $Rep_{G_\mathcal{Y}(C)}$. Because $G_\mathcal{Y}=\mathbb{G}_m$, this
extension is trivial in the category $Rep_{G_\mathcal{Y}}$. By the Tannakian
equivalence of categories, the extension $\mathcal{U} /\mathcal{W}$
is also trivial in $Ext_{\bold{T}}(\bold{1}, \mathcal{Y} /
\mathcal{W})$, and $\mathcal{W} \in \bold{V}$. Then  $\mathcal{V}
\subset \mathcal{W}$ by minimality. This concludes the
proof of Lemma \ref{lemme:minimal}, hence of Theorem \ref{theorem:Ber2}.

\begin{coro}\label{theorem:ind}
 Let
  $\mathcal{Y}$ be an object of $\bold{T}$,  let $\Delta$
  be the ring   $End(\mathcal{Y})$, and let
  $ \mathcal{E}_1,...,\mathcal{E}_n$ be  extensions  of $\bold{1}$
  by $ \mathcal{Y} $   such that
  $\mathcal{E}_1,...,\mathcal{E}_n$ are
  $\Delta$-linearly independent in
  $Ext^1_{\bold{T}}(\bold{1},\mathcal{Y})$. Assume that  $G_{\mathcal{E}_1}, ..., G_{\mathcal{E}_n}$ are reduced,
  that $G_\mathcal{Y} = \bold{G}_m$, and that the action of ${\bold{G}_m}$ on $\omega(\mathcal{Y})$ is given by its canonical character.
  Then the unipotent
radical of $G_{\mathcal{E}_1\oplus ... \oplus \mathcal{E}_n}$ is  isomorphic to $ \omega(\mathcal{Y})^n$.\\
\end{coro}

\textbf{Proof}

\medskip

For any extension $\mathcal{E}$ of  $\bold{1}$ by $ \mathcal{Y}$,
and for any $\alpha \in \Delta$, we denote by $\alpha_*
(\mathcal{E})$ the  pushout 
 of  $\mathcal{E}$ by $\alpha$ ;
this is how the structure of $\Delta$-module of
$Ext^1_{\bold{T}}(\bold{1}, \mathcal{Y})$ is defined.

\medskip

We first note that the direct sum $\mathcal{Y}^n$ admits  
$G_{\mathcal{Y}^n}= G_\mathcal{Y} = \mathbb{G}_m$ 
as Galois group, and that $\mathbb{G}_m$ again acts on $\omega(\mathcal{Y}^n) = \omega(\mathcal{Y})^n$ through its canonical character. On the other hand, the extension $\mathcal{E}_1 \oplus ...\oplus \mathcal{E}_n$ of $\bold{1}^n$ by $\mathcal{Y}^n$  and its pull-back $\mathcal{E} \in Ext^1_{\bold{T}}(\bold{1}, \mathcal{Y}^n)$ by the diagonal map from $\bold{1}$ to $\bold{1}^n$ generate in  $\bold T$ the same sub-Tannakian
 category. Therefore, their Galois groups  $G_{\mathcal{E}_1 \oplus... \oplus
\mathcal{E}_n}$ and $G_{\mathcal{E}}$ are equal, and reduced in view of our hypothesis. Let us assume that
the unipotent radical $R_u$ of $G_{\mathcal{E}}$ do not fill up
$ \omega(\mathcal{Y}^n) = \omega(\mathcal{Y})^n$.\\

By Theorem \ref{theorem:Ber2}, $R_u$  is equal to the $C$-vectorial group
$\omega(\mathcal{V})$ where $\mathcal{V} \in \bold{T}$ is the
smallest sub-object of $\mathcal{Y}^n$  such that  the quotient by
$\mathcal{V}$ of the extension $\mathcal{E}$ of $\bold{1}$ by
$\mathcal{Y}^n$ is trivial in the category $\bold{T}$. If
$\mathcal{V}$ is not equal to $\mathcal{Y}^n$, then
$\omega(\mathcal{V}) \subsetneq \omega(\mathcal{Y}^n)$. Because
$\omega(\mathcal{V})$ is a sub-representation of the representation
$\omega(\mathcal{Y}^n)$ of $G_\mathcal{Y} = \mathbb{G}_m$, it lies  in the kernel $H$
of a
 non trivial $G_\mathcal{Y}$-equivariant 
 homomorphism $\phi$ from  $\omega(\mathcal{Y}^n)$ to
 $\omega(\mathcal{Y})$. By tannakian equivalence of category,
 there then exists a non trivial morphism  $\Phi \in Hom_{\bold{T}}(\mathcal{Y}^n, \mathcal{Y})$
   such that $\mathcal{V} \subset \rm{Ker}(\Phi)$.\\  Now, consider the following diagram:
$$\xymatrix{
 \mathcal{Y}^n \ar[d] \ar[dr]^{\Phi}  \\
{\mathcal{Y}}^n / \mathcal{V} \ar[r] & {\mathcal{Y}}^n /
\rm{Ker}(\Phi)\simeq \mathcal{Y}.} $$ Since $\Phi \in Hom_{\bold{T}}(\mathcal{Y}^n, \mathcal{Y})$,
we can write $\Phi(X_1,...,X_n)= \alpha_1 X_1+...+\alpha_n X_n$,
with
  $\alpha_i \in End_{\bold{T}}(\mathcal{Y})$. Then $\Phi_{*}(\mathcal{E})={\alpha_1}_*
  (\mathcal{E}_1)
+{\alpha_2}_*(\mathcal{E}_2) +...+ {\alpha_n}_* (\mathcal{E}_n)$
is a quotient of $\mathcal{E}/ \mathcal{V}$, hence  a trivial  extension of $\bold{1}$ by $\mathcal{Y}$ in $\bold{T}$. In conclusion, the extension $\alpha_1 \mathcal{E}_1+...+\alpha_n
\mathcal{E}_n \in Ext^1_{\bold{T}}(\bold{1},\mathcal{Y})$ is
trivial. But this contradicts the $\Delta$-linearly independence  in
$Ext^1_{\bold{T}}(\bold{1},\mathcal{Y}) $  of the extensions $\mathcal{E}_1,..., \mathcal{E}_n$.

\section{Application to Theorem \ref{theorem:ilc3}}


 We shall apply Corollary \ref{theorem:ind} to the category $\bold{T} = \mathcal{T}$ of $t$-motives, which satisfies $C:=End_{\mathcal{T}}(\bold{1})   =\mathbb{F}_q(t)$, and all of whose Galois groups are reduced,  and to the Carlitz motive $\mathcal{Y} := {\cal C}$, for which  $\Delta :=  End_{\mathcal{T}}({\cal C})= \mathbb{F}_q(t)$ and $G_{\cal C} = \bold{G}_m$ acts on the line $\omega( {\cal C})$ through its canonical character. We recall the extensions $X(\alpha_i),  i = 1, ..., r$, of $\bold 1$  by $\cal C$ described in Section 1. Because of Corollary \ref{theorem:ind}, the dimension of the
algebraic group  $G=  G_{\bigoplus_{i=1}^rX(\alpha_i)}$
 is equal to $1
+ n$, where $n$ denotes the dimension of the
vector space over $\Delta=\mathbb{F}_q(t)$ generated by the $X(\alpha_i)$'s in $Ext^1_{\mathcal{T}}(\bold{1},{\cal C})$.\\

By an easy computation (similar to \cite{Hardouin}, 3.8), we get
$$n=max \lbrace s |\nexists f \in
\overline{k}(t), (\mu_i)_{i=1}^s \in \mathbb{F}_q(t) \mbox{ not all
zero}, \ \mbox{such that} (t-\theta)f^{(-1)}-f= \sum_{i=1}^s \mu_i
\alpha_i^{(-1)}(t-\theta)\rbrace.$$ By assumption,
$\tilde{\pi},log_{\cal C}(\alpha_1),...,log_{\cal C}(\alpha_r)$ are linearly
independent over $k=\mathbb{F}_q(\theta)$. Following \cite{Papa}, bottom of p. 50, we will now prove that under this hypothesis,  $n$ is equal
to $r$.\\

 Suppose that $n<r$. Then, let us consider $s$ such that $\exists f \in \overline{k}(t),
(\mu_i)_{i=1}^s \in \mathbb{F}_q(t)$  non all equal to zero such that 
\begin{equation}\label{eqn:ld}(t-\theta)f^{(-1)}-f= \sum_{i=1}^s \mu_i
\alpha_i^{(-1)}(t-\theta).\end{equation} It follows from Equation (\ref{eqn:ld}) that $f$ is regular at $t=\theta$ : if not,
$f^{(-1)}$ must have a pole at $t=\theta^{(-1)}$ which implies that
$f$ has a pole at $t=\theta^{(-1)}$. By repeating this argument, we
get  that if $f$ is singular at $t=\theta$ it is also singular at
$t=\theta^{(-i)}$ for all $i \geq 1$, which is impossible.
Therefore, $f$ and $f^{(-1)}$ are regular at $t=\theta$.\\

Considering the form of Equation (\ref{eqn:ld}), we then get
$f(\theta)=0$. Moreover, the solutions $y$ of  (\ref{eqn:ld}) are of the
following type :
$$ y =\mu \frac{1}{\Omega} + \sum_{i=1}^s \mu_i L_{\alpha_i}$$
with $\mu \in \mathbb{F}_q(t)$.
 So, there exists $\mu \in \mathbb{F}_q(t)$, such that :
\begin{equation}\label{eqn:ev}
f=\mu \frac{1}{\Omega} + \sum_{i=1}^s \mu_i L_{\alpha_i}
.\end{equation} By taking $t=\theta$ in (\ref{eqn:ev}), we get:
$$0=\mu(\theta)\tilde{\pi} + \sum_{i=1}^s \mu_i(\theta) log_{\cal C}(\alpha_i).$$
This  is a non trivial relation over $k$ between $ \tilde{\pi},
log_{\cal C}(\alpha_1),...,log_{\cal C}(\alpha_r)$, which contradicts our
assumption.\\

So, $dim
G=r+1$. This 
 concludes the
proof of Theorem \ref{theorem:ilc3}, and implies, as recalled in Section 1, that  $trdeg_{\overline{k}}\overline{k}(\tilde{\pi},log_C(\alpha_1),...,log_C(\alpha_r))= r+1$, i.e. that $\tilde{\pi},log_C(\alpha_1),...,log_C(\alpha_r)$ are
algebraically independent over $\overline{k}$.

\end{document}